\documentclass[12pt]{article}

\usepackage{amsmath} 
\usepackage{amssymb}
\usepackage{amsthm}
\usepackage{enumerate}

\newtheorem{thm}{Theorem}
\newtheorem{defn}{Definition}
\newtheorem{prop}[thm]{Proposition}
\newtheorem{lem}[thm]{Lemma}

\newcommand{\h}{\ensuremath{{\mathfrak h}}}

\newcommand{\Q}{\ensuremath{{\mathbb Q}}}
\newcommand{\Z}{\ensuremath{{\mathbb Z}}}
\newcommand{\Zeta}{\ensuremath{{\mathcal Z}}}

\allowdisplaybreaks[2]

\title{New Approach to Ohno Relation for Multiple Zeta Values\\
\small{ Dedicated to Professor Yoshiyuki Shimizu
on the occasion of his sixtieth birthday.}
}
\author{Jun-ichi OKUDA \and Kimio UENO}
\date{}

\usepackage{hyperref}
\begin{document}
\maketitle


\section{Introduction.\label{sec.intro}}
The multiple zeta values (MZVs) are real values
defined by 
\begin{equation*}
 \zeta(k_1,\dots,k_m)
 =
 \sum_{n_1 > n_2 > \dots >n_m>0} \frac{1}{n_1^{k_1}n_2^{k_2}\dotsb n_m^{k_m}}
  \quad \left(k_1,\dots,k_m \in \Z_{\ge 1}, k_1 \ge 2\right).
\end{equation*}
The weight and the depth of $\zeta(k_1,\dots,k_m)$ are
$k_1 + \dotsb + k_m$ and $m$, respectively.

Recently, MZVs have been studied extensively
in number theory \cite{gonc,zagi1},
knot theory \cite{le.mura},
mirror symmetry \cite{hoff3}
and perturbative quantum field theory \cite{krei}.
Many relations of MZVs,
for example
Hoffman's relation \cite{hoff1},
the duality formula \cite{zagi1},
the sum formula \cite{gran},
Le-Murakami's relation \cite{le.mura}
and the cyclic sum formula \cite{hoff.ohno}, have been discovered.
In particular,
Ohno discovered a systematic relation
which contain Hoffman's relation, the duality formula and the sum formula.
These relations are important to consider the structure of \Zeta,
which is an algebra over $\Q$ generated by all the MZVs.

\begin{thm}[Ohno's relation~\cite{ohno}]\label{thm.ohno}
 Any index ${\mathbf k} = (k_1,\dots,k_m) \in \Z_{\ge 1}^m$, $k_1\ge 2$,
 can be written for some $s\in\Z_{\ge 1}$
 and $a_i,b_i\in\Z_{\ge 1}\,(i=1,\dotsc,s)$ as 
 \begin{align*}
  {\mathbf k}  &= 
  (a_1 + 1, \underbrace{1,\dots,1}_{b_1 - 1},
   a_2 + 1, \underbrace{1,\dots,1}_{b_2 - 1},
   \dots,
    a_s + 1, \underbrace{1,\dots,1}_{b_s - 1}).\\
   \intertext{ Let ${\mathbf k'} = (k_1',\dots,k_{m'}')$ be the dual index
                for ${\mathbf k}$ defined by}\\
  {\mathbf k}'  &= 
  (b_s + 1, \underbrace{1,\dots,1}_{a_s - 1},
   b_{n-1} + 1, \underbrace{1,\dots,1}_{a_{n-1} - 1},
   \dots,
    b_1 + 1, \underbrace{1,\dots,1}_{a_1 - 1}).
 \end{align*}
Then for all $l\in\Z_{\ge 0}$,
 \begin{multline*}
  \sum_{\substack{
        \epsilon_1,\dots,\epsilon_m\ge 0\\
        \epsilon_1 + \dots + \epsilon_m = l}}
   \!\!
   \zeta(k_1+\epsilon_1,k_2+\epsilon_2,\dots,k_m+\epsilon_m)\\
  \qquad\qquad\qquad=
   \!\!
  \sum_{\substack{
        \epsilon_1',\dots,\epsilon_{m'}'\ge 0\\
        \epsilon_1' + \dots + \epsilon_{m'}' = l}}
   \!\!
   \zeta(k_1'+\epsilon_1',k_2'+\epsilon_2',\dots,k_{m'}'+\epsilon_{m'}')
 \end{multline*}
\end{thm}

Ohno's proof is due to iterated integral expression of
the generating functions in $X_1,\dotsc,X_s$
for given $l\in\Z{\ge 0}$ and ${\mathbf k}$
\begin{multline*}
 \sum_{\substack{l_i\ge0\\\l_1+\dotsb+l_s=l}}\ 
  \idotsint\limits_{1>t_1>\dotsb>0}
   \prod_{i=1}^s
    \left(\frac{dt}{t}\right)^{a_i}
    \left(\frac{dt}{t} + X_i \frac{dt}{1-t}\right)^{b_i+l_i-1}\!
    \frac{dt}{1-t}\\
 =
  \sum_{\substack{l_i\ge0\\\l_1+\dotsb+l_s=l}}
  \sum_{1\le d_i \le b_i + l_i}
   S_{\mathbf k}(l_1,\dotsc,l_s;d_1,\dotsc,d_s)\
    X_1^{d_1-1}\dotsb X_s^{d_s-1}.
\end{multline*}
The coefficient of $X_1^{b_1}\dotsb X_s^{b_s}$ is
the left hand of Theorem~\ref{thm.ohno}.
%where the coefficient $S({\mathbf k};d_1,\dotsc,d_s)$ is defined by
%\begin{align*}
% S({\mathbf k};d_1,\dotsc,d_s)
% =&
%  \!\!\!\!  \!\!\!\!
% \sum_{\substack{\delta_{i,j}\ge 1\\\sum \delta_{i,j}= l + \sum b_i}}
%  \!\!\!\!  \!\!\!\!
%  \zeta(
%   \underbrace{a_1+\delta_{1,1},\delta_{1,2},\dotsc,\delta_{1,d_1}}_{d_1},
%    \dots,
%     \underbrace{a_s+\delta_{s,1},\delta_{s,2},\dotsc,\delta_{s,d_s}}_{d_s}).
%\end{align*}

We consider generating functions
of left hand side of Theorem~\ref{thm.ohno}
with respect to the increase $l$ of weight
\begin{equation*}
 \sum_{l=0}^\infty
  \left\{
   \sum_{\substack{
        \epsilon_1,\dots,\epsilon_m\ge 0\\
        \epsilon_1 + \dots + \epsilon_m = l}}
   \!\!\zeta(k_1+\epsilon_1,k_2+\epsilon_2,\dots,k_m+\epsilon_m)
  \right\} \lambda^l.
\end{equation*}
We shall find functional relations among them.
Using these relations,
we give a new proof of Ohno's relation.

This article is organized as follows:
In Section~\ref{sec.harm}, we restate Ohno's relation
and introduce our generating function via the harmonic algebra.
In Section~\ref{sec.rel}, we state and prove
the functional relations of generating functions.
In Section~\ref{sec.proof}, we prove Ohno's relation.

\section{A translation of Ohno's relation. \label{sec.harm}}

We explain Ohno's relation
in terms of the derivations~\cite{kane1, hoff.ohno}
introduced to the harmonic algebra~\cite{hoff2}.

\begin{defn}
Let $\h = \Q\langle x,y\rangle$
be a noncommutative polynomial algebra over $\Q$
generated by $x$ and $y$.
For any monomial $w$ in $\h$, weight and depth mean
the length of w and the number of y which appears in w, respectively.

We define a new multiplication $*$ of $\h$ as follows{\normalfont :}
For any words $w$,$w_1$,$w_2$ $\in \h$
and positive integers $p,q \in \Z_{\ge 1}$,
we set
  \begin{enumerate}[{\normalfont\rmfamily (H.1)}]
  \item $1*w = w*1 = w$
  \item $x^p*w = w*x^p = wx^p$
  \item $x^pyw_1*x^qyw_2 =
             x^py(w_1*x^qyw_2) + x^{(p+q+1)}y(w_1*w_2)
             +x^qy(x^pyw_1*w_2)$
  \end{enumerate}
We call the triple $(\h, +, *)$ harmonic algebra.
\end{defn}

Hoffman showed the next theorem which relates the harmonic
algebra with the algebra $\Zeta$.

\begin{thm}[Hoffman~\cite{hoff2}]
The harmonic algebra $(\h, +, *)$ is a
commutative polynomial algebra of infinite indeterminates,
and for a subalgebra $\h^0 := \Q + x\,\h\,y$ of $\h$,
there exists an algebra homomorphism
\begin{equation*}
	   \widetilde{\zeta}: \h^0 \longrightarrow \Zeta
\end{equation*} 
such that
\begin{equation*}
  \widetilde{\zeta}(x^{k_1}y^{l_1}\dots x^{k_m}y^{l_m}) = 
  \zeta(k_1 + 1,\underbrace{1,\dots,1}_{l_1-1},
   \dots,
   k_m+1,\underbrace{1,\dots,1}_{l_m-1}).
\end{equation*}
\end{thm}

By virtue of this theorem, relations of MZVs are thought of to be the
images of elements in $\ker \widetilde{\zeta}$.
In this context,
Ohno's relation is restated
in terms of the image of certain derivations on $\h$.
To describe these derivations,
we introduce some notations.

\begin{defn}[\cite{hoff2,kane1,hoff.ohno}]\label{defn.map}
 Let $\tau$ be an anti-involution on $\h$ defined by
 \begin{align*}
  \tau(x) = y,\qquad \tau(y) = x.
 \end{align*}
 For for any $n \in \Z_{\ge 1}$,
 let $D_n$ be a derivation of $\h$ defined by
 \begin{align*}
  D_n(x) = 0,\qquad D_n(y)= x^n y.
 \end{align*}
\end{defn}

Through the homomorphism $\widetilde{\zeta}$,
 the action of $D_n$ is regarded as follows:
\begin{align*}
 \widetilde{\zeta}\left(D_n(x^{k_1-1}y\dotsm x^{k_m-1}y)\right)
  = \sum_{i=1}^m \zeta(k_1,\dots,k_i + n,\dots,k_m).
\end{align*}

\begin{prop}[Ohno's relation \cite{kane2,hoff.ohno}]
Let $\lambda$ be a parameter.\\
Then in $\h[[\lambda]]$, a ring of formal power series in $\lambda$
over $\h$, we have
 \begin{align*}
  \text{coefficients of } 
   \left(
    \exp\left(\sum_{n=0}^\infty \frac{D_n}{n}\lambda^n\right) 
     -\exp\left(\sum_{n=0}^\infty \frac{D_n}{n}\lambda^n\right)\circ\tau
   \right)(\h^{0})
   \subset\ker\widetilde{\zeta}.
 \end{align*}
\end{prop}

As $\sum_{n=0}^\infty \frac{D_n}{n}\lambda^n$ is
a derivation on $\h[[\lambda]]$,
$\exp\left(\sum_{n=0}^\infty \frac{D_n}{n}\lambda^n\right)$ is
an automorphism on $\h[[\lambda]]$.
So we can easily calculate this version of Ohno's relation.
By definition of $D_n$,
$\exp\left(\sum_{n=0}^\infty \frac{D_n}{n}\lambda^n\right)$ sends $x$ and $y$
to $x$ and $(1-x\lambda)^{-1}y$ respectively.
Set $L_i := l_1 + \dots + l_i (i\ge 1)$ and $L_0 := 0$ for short, and
\begin{align*}
 &\left(\exp(\sum_{n=0}^\infty \frac{D_n}{n}\lambda^n)\right)
   (x^{k_1}y^{l_1}x^{k_2}y^{l_2}\dotsm x^{k_m}y^{l_m})\\
 &\quad
  = x^{k_1} ((1-x\lambda)^{-1} y)^{l_1}
    x^{k_2} ((1-x\lambda)^{-1} y)^{l_2}
     \dotsm
    x^{k_m} ((1-x\lambda)^{-1} y)^{l_m}\\
 &\quad
  = x^{k_1} \Bigl(\sum_{j=0}^\infty x^jy\lambda^j\Bigr)^{l_1}
    x^{k_2} \Bigl(\sum_{j=0}^\infty x^jy\lambda^j\Bigr)^{l_2}
     \dotsm
    x^{k_m} \Bigl(\sum_{j=0}^\infty x^jy\lambda^j\Bigr)^{l_m}\\
 &\quad
  = 
  \sum_{j_1=0}^\infty \dotsm \!\!\!\sum_{j_{L_m}=0}^\infty \!\!\!
  \begin{aligned}[t]
  &x^{k_1}\left(x^{j_1}y\lambda^{j_1}\dotsm x^{j_{L_1}}y\lambda^{j_{L_1}}\right)
   \cdot
    x^{k_2}\left(x^{j_{L_1+1}}y\lambda^{j_{L_1+1}}
     \dotsm x^{j_{L_2}}y\lambda^{j_{L_2}}\right)\\
  &\qquad\qquad\dotsm
   x^{k_m}\left(x^{j_{L_{m-1}+1}}y\lambda^{j_{L_{m-1}+1}}
                 \dotsm x^{j_{L_m}}y\lambda^{j_{L_m}}\right)
 \end{aligned}\\
 &\quad
 = \sum_{n=0}^\infty
  \left\{
   \sum_{j_1+\dots+j_{L_m}=n}
 \begin{aligned}[t]
    &x^{k_1}(x^{j_1}y\dotsm x^{j_{L_1}}y)
    \cdot x^{k_2}(x^{j_{L_1+1}}y\dotsm x^{j_{L_2}}y)\\
    &\qquad\qquad\qquad\dotsm
    x^{k_m}(x^{j_{L_{m-1}+1}}y\dotsm x^{j_{L_m}}y)
 \end{aligned}  
\right\}
  \lambda^n.
\end{align*}
Hence we have
\begin{align}
  &\left(\widetilde{\zeta}\circ\exp\left(
       \sum_{n=0}^\infty \frac{D_n}{n}\lambda^n
      \right)\right)
   (x^{k_1}\underbrace{y\dotsm y}_{l_1}
    \dotsm
    x^{k_m}\underbrace{y\dotsm y}_{l_m})\notag\\
 &\quad=
 \begin{aligned}[t]
 &\sum_{n=0}^\infty
 \Biggl\{
  \sum_{j_1+\dotsb+j_{L_m}=n}
    \zeta(\underbrace{k_1+1+j_1,1+j_2,\dotsc,1+j_{L_1}}_{l_1},\\
    &\qquad\qquad\qquad\dotsc,
     \underbrace{k_m+1+j_{L_{m-1}+1},1+j_{L_{m-1}+2},\dotsc,1+j_{L_m}}_{l_m})
 \Biggr\} \lambda^n
   \end{aligned}\notag\\
 \begin{split}
 &\quad=
  \sum_{n_1 > \dots > n_{L_m > 0}}
  \frac{1}{n_1^{k_1}}
  \underbrace{
   \sum_{j_1=1}^\infty\frac{\lambda^{j_1-1}}{n_1^{j_1}}
   \dots
   \sum_{j_{l_1}=1}^\infty\frac{\lambda^{j_{l_1}-1}}{n_{l_1}^{j_{l_1}}}
  }_{l_1}\\
  &\qquad\qquad\qquad\times\dotsm\times
  \frac{1}{n_{L_{m-1}+1}^{k_m}}
  \underbrace{
   \sum_{j_{L_{m-1}+1}=1}^\infty
    \frac{\lambda^{j_{L_{m-1}+1}-1}}%
     {n_{L_{m-1}+1}^{j_{L_{m-1}+1}}}
   \dots
   \sum_{j_{L_m}=1}^\infty
    \frac{\lambda^{j_{L_m}-1}}{n_{L_m}^{j_{L_m}}}
  }_{l_m}
 \end{split} \tag{$*$} \label{power.series}\\
\begin{split}
 &\quad=
    \sum_{n_1 > \dots > n_{L_m > 0}}
   \frac{1}{
    n_1^{k_1}
     \underbrace{ (n_1-\lambda)\dotsm (n_{L_1}-\lambda)}_{l_1}
    n_{L_1+1}^{k_2}
     \underbrace{ (n_{L_1+1}-\lambda)\dotsm (n_{L_2}-\lambda)}_{l_2}}\\
   &\qquad\qquad\qquad\times\dotsm\times
   \frac{1}{
    n_{L_{m-1}+1}^{k_m}
     \underbrace{ (n_{L_{m-1}+1}-\lambda)\dotsm(n_{L_m}-\lambda)}_{l_m}
  }. 
\end{split}
\tag{$**$} \label{partial.series}
\end{align}
The power series \eqref{power.series} absolutely converges for $|\lambda| < 1$.
The series \eqref{partial.series} absolutely converges
at points except positive integers.

%and defines a meromorphic function on $\C$
%with poles at positive integers.
%
\begin{defn}
For $m \in \Z_{\ge 1}$,
let us consider a sequence
$\{k_i,l_i\}_{i=1}^m \in \Z_{\ge 0}^{2m}$, $k_1, l_m \ge 1$.
We set $L_i = l_1 + \dotsb + l_i$ and $K_i = k_m + \dotsb + k_i$
and define generating functions of MZVs as follows{\normalfont :}
\begin{align*}
 f(\{k_i,l_i\}_{i=1}^m;\lambda)
 :=& \left(
      \widetilde{\zeta}
       \circ \exp\left(\sum_{n=0}^\infty \frac{D_n}{n}\lambda^n \right)
     \right)
      \left(x^{k_1}y^{l_1}\dotsm x^{k_m}y^{l_m}\right)\\
 =& 
   \sum_{n_1 > \dots > n_{L_m > 0}}
   \prod_{i=1}^m
    \frac{1}{n_{L_{i-1}+1}^{k_i}
     \underbrace{(n_{L_{i-1}+1}-\lambda)\dotsm(n_{L_i}-\lambda)}_{l_i}
    },
  \\
 g(\{k_i,l_i\}_{i=1}^m;\lambda)
 :=& \left(
      \widetilde{\zeta}
       \circ \exp\left(\sum_{n=0}^\infty \frac{D_n}{n}\lambda^n \right)
     \right)
      \left(\tau(x^{k_1}y^{l_1}\dotsm x^{k_m}y^{l_m})\right)\\
 =&
   \sum_{0< n_1 < \dots < n_{L_m}}
   \prod_{i=m}^1
    \frac{1}{n_{K_{i-1}+1}^{l_i}
     \underbrace{(n_{K_i}-\lambda)\dotsm(n_{K_{i-1}+1}-\lambda)}_{k_i}.
    }
\end{align*}
It is obvious from the definition of $\tau$ that
\begin{align*}
 g(\{k_i,l_i\}_{i=1}^m;\lambda) = f(\{l_i,k_i\}_{i=m}^1;\lambda).
\end{align*}
\end{defn}
We show that these generating functions have 
partial-fractions expansion with simple poles at positive integers.

\begin{lem}\label{key.lem}
The generating function $f(\{k_i,l_i\}_{i=1}^m;\lambda)$
can be written as follows{\normalfont :}
 \begin{align*}
f(\{k_i,l_i\}_{i=1}^m;\lambda)=
  \sum_{n=1}^\infty
   \left\{
    \sum_{j=1}^{L_m}
    \sum_{\substack{n_1 > \dots > n_{j-1} > n\\ n >n_{j+1}>\dots>n_{L_m}}}
     C_n^{n_1\dotsc \stackrel{\stackrel{j}\smallsmile}n \dotsc n_{L_m}}\right\}
    \frac{1}{n-\lambda}
\end{align*}
where
\begin{align*}
C_{n_j}^{n_1\dotsc n_{L_m}}
  = \frac{1}{n_1^{k_1}n_{L_1+1}^{k_2}\dotsc n_{n_{L_{m-1}}+1}^{k_m}}
          \prod_{i\not=j}\frac{1}{(n_i-n_j)}.
 \end{align*}
\end{lem}

\begin{proof}
The generating function $f$ can be written as follows:
\begin{align*}
f(\{k_i,l_i\};\lambda)
 = 
  \sum_{n_1 > \dots > n_{L_m > 0}}
  \sum_{j=1}^{L_m}
   \frac{C_{n_j}^{n_1\dotsc n_{L_m}}}{n_j-\lambda}.
\end{align*}
For the proof, we have to show that it is possible to
change the order of the summation.
So it is sufficient to prove that for any $j$
\begin{align*}
   \sum_{n_1 > \dots > n_{L_m > 0}}
   \frac{C_{n_j}^{n_1\dotsc n_{L_m}}}{n_j-\lambda}
\end{align*}
converges absolutely.
It can be proved by virtue of the idea of Mordell~\cite{mord}.
Put $d_i = n_i - n_{i+1}$ for $i=1,\dotsc,L_m-1$ and $d_{L_m}= n_{L_m}$.
Making use of the inequality
\begin{align*}
 d_1 + d_2 + \dotsb + d_{L_m} \ge L_m \sqrt[L_m]{d_1d_2\dotsb d_{L_m}}
\end{align*}
we have
\begin{align*}
 \left|\frac{C_{n_j}^{n_1\dotsc n_{L_m}}}{n_j-\lambda}\right|
 &=\left|
   \left\{
     \frac{1}{n_1^{k_1}\dotsm n_{n_{L_{m-1}}+1}^{k_m}}
          \prod_{i\not=j}\frac{1}{(n_i-n_j)}\right\}
    \frac{1}{n_j-\lambda}
  \right|\\
 &=
 \frac{1}{(d_1+\dotsb+d_{L_m})^{k_1}
           \dotsb(d_{L_{m-1}+1}+\dotsb+d_{L_m})^{k_m}}\\
 &\qquad\qquad\times
  \left\{
   \prod_{i<j}\frac{1}{(d_i+\dotsb+d_{j-1})}
   \prod_{i>j}\frac{1}{(d_j+\dotsb+d_{i-1})}
  \right\}\\
 &\qquad\qquad\qquad\qquad\times
  \frac{1}{\left|d_j+\dotsb+d_{L_m}-\lambda\right|}\\
 &\le
  \frac{1}{\left(L_m \sqrt[L_m]{d_1\dotsb d_{L_m}}\right)^{k_1}}
   \left\{\prod_{i=1}^{L_m-1}\frac{1}{d_i}\right\}
    \frac{1}{\left|d_j+\dotsb+d_{L_m}-\lambda\right|}
\end{align*}
Let $\lambda$ be in a compact set which dose not involve positive
 integers.
Then there exists a positive constant $A$ such that
\begin{align*}
 \frac{1}{\left|d_j+\dotsb+d_{L_m}-\lambda\right|}
 \le
 \frac{A}{d_{L_m}}.
\end{align*}
Hence
\begin{align*}
 &\sum_{n_1 > \dots > n_{L_m > 0}}
  \left|\frac{C_{n_j}^{n_1\dotsc n_{L_m}}}{n_j-\lambda}\right|\\
 &\qquad\qquad\le
 A
 \sum_{d_1,\dotsc,d_{L_m}=1}^\infty
  \frac{1}{\left(L_m \sqrt[L_m]{d_1\dots d_{L_m}}\right)^{k_1}}
   \left\{\prod_{i=1}^{L_m}\frac{1}{d_i}\right\}\\
 &\qquad\qquad<+\infty.
\end{align*}

\end{proof}

Using $f$ and $g$, Ohno's relation is restated as follows.

\begin{thm}[Ohno's relation, again]
 For any sequence
 $\{k_i,l_i\}_{i=1}^m$$\in$$\Z_{\ge 0}^{2m}$ with  $k_1$,$l_m$$\ge$$1$,
 \begin{equation*}
  f(\{k_i,l_i\}_{i=1}^m;\lambda) = g(\{k_i,l_i\}_{i=1}^m;\lambda).
 \end{equation*}
\end{thm}

\section{Functional relations of $f$'s.\label{sec.rel}}

\begin{thm}\label{mainTh}
 We set
 $\lambda' := \lambda - 1$ and 
 $I:=\left\{(0,0),(1,0),(0,1)\right\}$,
 then the generating function $f$ satisfies the following relations.
 \begin{enumerate}[{\normalfont \upshape (i)}]
 \item \label{thm.1} If $k_1,\,l_m \not = 1$,
  \begin{align*}
   &\sum_{\{(\delta_i,\epsilon_i)\}_{i=1}^m\in I^m}
    (-\lambda)^{m-|\delta|-|\epsilon|}
    f(\{k_i-\delta_i,l_i-\epsilon_i\}_{i=1}^m;\lambda)\\
   &= 
   \sum_{\delta_1'\in\{0,1\}}
   \sum_{\{(\delta_j',\epsilon_j')\}_{j=2}^m\in I^{m-1}}
   \sum_{\epsilon_{m+1}'\in\{0,1\}}
    (-\lambda')^{m-|\delta'|-|\epsilon'|}
     f(\{k_i-\delta_i',l_i-\epsilon_{i+1}'\}_{i=1}^m;\lambda').
  \end{align*}
 \item If $k_1 =1,\,l_m \not = 1$,
  \begin{align*}
   &\sum_{\epsilon_1\in\{0,1\}}
    \sum_{\{(\delta_i,\epsilon_i)\}_{i=2}^m\in I^{m-1}}
    (-\lambda)^{m-|\delta|-|\epsilon|}
    f(\{1,l_1-\epsilon_1\}\cup
          \{k_i-\delta_i,l_i-\epsilon_i\}_{i=2}^m;\lambda)\\
   &=
   \sum_{\{(\delta_j',\epsilon_j')\}_{j=2}^m\in I^{m-1}}
   \sum_{\epsilon_{m+1}'\in\{0,1\}}
   \begin{aligned}[t]
    &(-\lambda')^{m-|\delta'|-|\epsilon'|}\\
    &\quad\times
     f(\{1,l_1-\epsilon_2'\}\cup\{k_i-\delta_i',l_i-\epsilon_{i+1}'\}_{i=2}^m
       ;\lambda').
   \end{aligned}  
\end{align*}
 \item If $k_1 \not = 1,\,l_m = 1$,
  \begin{align*}
   &\sum_{\{(\delta_i,\epsilon_i)\}_{i=1}^{m-1}\in I^{m-1}}
    \sum_{\delta_m\in\{0,1\}}
    (-\lambda)^{m-|\delta|-|\epsilon|}
    f(\{k_i-\delta_i,l_i-\epsilon_i\}_{i=1}^{m-1}\cup\{k_m-\delta_m,1\}
       ;\lambda)\\
   &= 
   \sum_{\delta_1'\in\{0,1\}}
   \sum_{\{(\delta_j',\epsilon_j')\}_{j=2}^m\in I^{m-1}}
   \begin{aligned}[t]
     &(-\lambda')^{m-|\delta'|-|\epsilon'|}\\
     &\quad\times
    f(\{k_i-\delta_i',l_i-\epsilon_{i+1}\}_{i=1}^{m-1}\cup\{k_m-\delta_m',1\}
       ;\lambda').
   \end{aligned}  
\end{align*}
 \item If $k_1 = l_m = 1$,
  \begin{align*}
   &\sum_{\epsilon_1\in\{0,1\}}
    \sum_{\{(\delta_i,\epsilon_i)\}_{i=2}^{m-1}\in I^{m-2}}
    \sum_{\delta_m\in\{0,1\}}
    (-\lambda)^{m-|\delta|-|\epsilon|}
    f(\{1,l_1-\epsilon_1,\dotsc,k_m-\delta_m,1\};\lambda)\\
   &= 
   \sum_{\{(\delta_j',\epsilon_j')\}_{j=2}^m\in I^{m-1}}
     (-\lambda')^{m-|\delta'|-|\epsilon'|}
    f(\{1,l_1-\epsilon_2',\dotsc,k_m-\delta_m',1\};\lambda').
  \end{align*}
 \end{enumerate}
 Here $|\delta^{(\prime)}|$ (resp. $|\epsilon^{(\prime)}|$) is
 the sum of all $\delta_i^{(\prime)}$ (resp. $\epsilon_i^{(\prime)}$).
 The generating function g also satisfies the same relations.
\end{thm}

To prove this theorem, we need a lemma.

\begin{defn}
 For $k_i,l_i\in\Z_{\ge 0}(i=1,\dots,m),\,k_1,l_m\ge1$ ,we set
 \begin{align*}
  \begin{split}
   &[\{(n-a_i)^{k_i},l_i\}_{i=1}^m;\lambda]\\
  &\qquad\qquad
   := \sum_{n_1 > \dots > n_{L_m} > 0}
   \prod_{i=1}^m
   \frac{1}{(n_{L_{i-1}+1}-a_i)^{k_i}
    \underbrace{
     (n_{L_{i-1}+1}-\lambda)\dotsm(n_{L_{i}}-\lambda)
    }_{l_i}
   },
  \end{split}   
 \end{align*}
where we interpret special cases with $k_i=0$ or $l_i=0$ for some $i$
as follows{\normalfont :}
 \begin{align*}
   [\{\ldots,l_{i-1},n^0,l_i,\ldots\};\lambda]
   &= [\{\ldots,l_{i-1}+l_i,\ldots\};\lambda], \\
   [\{\ldots,(n-a)^{k_{i-1}},0,(n-a)^{k_i},\ldots\};\lambda]
   &= [\{\ldots,(n-a)^{k_{i-1}+k_i},\ldots\};\lambda].
 \end{align*}
\end{defn}

Using this notation, $f$ and $g$ are expressed as follows:
\begin{align*}
 \left\{
 \begin{aligned}[c]
   &f(\{k_i,l_i\}_{i=1}^m;\lambda) = [\{n^{k_1},l_i\}_{i=1}^m;\lambda], \\
  &g(\{k_i,l_i\}_{i=1}^m;\lambda) = [\{n^{l_1},k_i\}_{i=m}^1;\lambda].
 \end{aligned}
 \right.
\end{align*} 

\begin{lem}\label{lemma}
\begin{enumerate}[{\normalfont \upshape (i)}]
\item If $m\not = 1$ or $l_m \not = 1$
\begin{enumerate}[{\normalfont \upshape (a)}]
\item\label{lemma.2} If $k_1\ge2$,
 \begin{multline*}
   \lambda\,[\{n^{k_1},l_1,\dots\};\lambda] 
   - [\{n^{k_1-1},l_1,\dots\};\lambda]
   - [\{n^{k_1},l_1-1,\dots\};\lambda]\\
   = 
   \lambda'\,[\{(n-1)^{k_1},l_1,\dots\};\lambda]
    - [\{(n-1)^{k_1-1},l_1,\dots\};\lambda].
 \end{multline*}
\item\label{lemma.3} If $k_1=1$,
 \begin{equation*}
  \lambda\,[\{n^1,l_1,\dotsc\};\lambda] - [\{n^1,l_1-1,\dotsc\};\lambda]
   = \lambda'\,[\{(n-1)^1,l_1,\dotsc\};\lambda].
 \end{equation*}
\end{enumerate}
\item\label{lemma.1} If $i\not=1$ and $i\not=m$ or $i=m$ and $l_m\not=1$,
  \begin{align*}
   &\lambda\,[\{\dots,n^{k_i},l_i,\dots\};\lambda] \\
   &\qquad
    - [\{\dots,n^{k_i-1},l_i,\dots\};\lambda]
    - [\{\dots,n^{k_i},l_i-1,\dots\};\lambda]\\
   &=
     \lambda'\,[\{\dots,l_{i-1},(n-1)^{k_i},\dots\};\lambda]\\
     &\qquad
      - [\{\dots,l_{i-1},(n-1)^{k_i-1},\dots\};\lambda]
      - [\{\dots,l_{i-1}-1,(n-1)^{k_i},\dots\};\lambda].
  \end{align*}
\item
\begin{enumerate}[{\normalfont \upshape (a)}]
\setcounter{enumii}{2}
\item\label{lemma.4} If $l_m=1$,
 \begin{align*}
   &\lambda\,[\{(n-1)^{k_1},l_1,(n-1)^{k_2},l_2,\dots,
      (n-1)^{k_{m-1}},l_{m-1},n^{k_m},1\};\lambda]\\
   &\qquad- [\{(n-1)^{k_1},l_1,(n-1)^{k_2},l_2,\dots,
      (n-1)^{k_{m-1}},l_{m-1},n^{k_m-1},1\};\lambda]\\
   &=
    \lambda'\,[\{n^{k_1},l_1,n^{k_2},l_2,\dots,n^{k_m},1\};\lambda']
    - [\{n^{k_1},l_1,n^{k_2},l_2,\dots,n^{k_m-1},1\};\lambda']\\
   & \qquad\qquad
     -[\{n^{k_1},l_1,n^{k_2},l_2,\dots,l_{m-1}-1,n^{k_m},1\};\lambda'].
 \end{align*} 
\item\label{lemma.5} If $l_m \ge 2$,
 \begin{multline*}
  [\{(n-1)^{k_1},l_1,(n-1)^{k_2},l_2,\dots,(n-1)^{k_m},l_m\};\lambda]\\
  = [\{n^{k_1},l_1,\dots,n^{k_m},l_m\};\lambda']
  - \frac{1}{\lambda'} [\{n^{k_1},l_1,\dots,n^{k_m},l_m-1\};\lambda'].
 \end{multline*}
\end{enumerate}
\end{enumerate}
\end{lem}

\begin{proof}
We use the partial-fractions expansion:
\begin{align}
 \begin{split}
  &\frac{\lambda}{n^k(n-\lambda)} - \frac{1}{n^{k-1}(n-\lambda)}\\
  &\qquad= 
  \frac{\lambda'}{(n-1)^k(n-\lambda)} - \frac{1}{(n-1)^{k-1}(n-\lambda)}
  + \left(\frac{1}{(n-1)^k} - \frac{1}{n^k} \right),
 \end{split}
 \tag{$\spadesuit$}\label{frac.1}\\
 &\frac{\lambda}{n(n-\lambda)}
 = \frac{\lambda'}{(n-1)(n-\lambda)}
 + \left(\frac{1}{n-1} - \frac{1}{n}\right).
 \tag{$\heartsuit$}\label{frac.2}
\end{align}

\begin{enumerate}[{\normalfont \upshape (i)}]
\item
\begin{enumerate}[{\normalfont \upshape (a)}]
\item
We set $B$ by
\begin{equation*}
 \prod_{j=1}^m   
   n_{L_{j-1}+1}^{k_j}
    \underbrace{
     (n_{L_{j-1}+1}-\lambda)\dotsm(n_{L_{j}}-\lambda)
    }_{l_j}
   = n_1^{k_1}(n_1-\lambda)\,B.
\end{equation*}
Then using \eqref{frac.1} we have
 \begin{align*}
  &\lambda\,[\{n^{k_1},l_1,\dotsc\};\lambda]
    - [\{n^{k_1-1},l_1,\dotsc\};\lambda]\\
  &\qquad=
   \sum_{n_1 > \dots>n_{L_m} > 0}
   \left\{
    \frac{\lambda}{n_1^{k_1}(n_1-\lambda)}
    - \frac{1}{n_1^{k_1-1}(n_1-\lambda)}
   \right\}\frac{1}{B}\\
  &\qquad=
   \sum_{n_1 > \dots>n_{L_m} > 0}
   \Biggl\{
    \frac{\lambda'}{(n_1-1)^{k_1}(n_1-\lambda)}
    - \frac{1}{(n_1-1)^{k_1-1}(n_1-\lambda)}\\
    & \hspace{8cm}
    + \left(\frac{1}{(n_1-1)^{k_1}}-\frac{1}{n_1^{k_1}}\right)
   \Biggr\}\frac{1}{B}\\
  \begin{split}
   &\qquad=
   \sum_{n_1 > \dots>n_{L_m} > 0}
   \left\{
    \frac{\lambda'}{(n_1-1)^{k_1}(n_1-\lambda)}
    - \frac{1}{(n_1-1)^{k_1-1}(n_1-\lambda)}
   \right\}\frac{1}{B}\\
   &\hspace{4cm}
    +\sum_{n_2>\dots>n_{L_m}>0}\,\sum_{n_1 = n_2+1}^\infty
     \left(\frac{1}{(n_1-1)^{k_1}}-\frac{1}{n_1^{k_1}}\right)\frac{1}{B}
  \end{split}\\
  &\qquad=
  \begin{aligned}[t]
   &\lambda'\,[\{(n-1)^{k_1},l_1,\dotsc\};\lambda]
   -[\{(n-1)^{k_1-1},l_1,\dotsc\};\lambda]\\
    &\hspace{7cm}
     +\sum_{n_2>\dots>n_{L_m}>0}
      \frac{1}{n_2^{k_1}B}
  \end{aligned}\\
  &\qquad=
  \begin{aligned}[t]
   &\lambda'\,[\{(n-1)^{k_1},l_1,\dotsc\};\lambda]
    -[\{(n-1)^{k_1-1},l_1,\dotsc\};\lambda]\\
   &\hspace{7cm}
    +[\{n^{k_1},l_1-1,\dotsc\};\lambda].
  \end{aligned}
 \end{align*}

\item
 Using \eqref{frac.2}, it can be proved in the same manner as \eqref{lemma.2}.
\end{enumerate}
\item We set $A$ and $B$ by
  \begin{align*}
   &A := \prod_{j=1}^{i-1}
   n_{L_{j-1}+1}^{k_j}
    \underbrace{
     (n_{L_{j-1}+1}-\lambda)\dotsm(n_{L_{j}}-\lambda)
    }_{l_j}, \\
   &B:=
   \Biggl(
    \prod_{j=i-1}^m
    n_{L_{j-1}+1}^{k_j}
     \underbrace{
      (n_{L_{j-1}+1}-\lambda)\dotsm(n_{L_{j}}-\lambda)
     }_{l_j}
   \Biggr) \Bigg/
    n_{L_{i-1}+1}^{k_i}(n_{L_{i-1}+1}-\lambda).
  \end{align*}
Then
 \begin{align*}
  &\lambda [\{\dots,l_{i-1},n^{k_i},l_i,\dots\};\lambda]
   - [\{\dots,l_{i-1},n^{k_i-1},l_i,\dots\};\lambda]\\
  &\quad=
   \begin{aligned}[t]
     &\sum_{n_1 > \dots > n_{L_m} >0}
     \frac{\lambda}{A\,n_{L_{i-1}+1}^{k_i}(n_{L_{i-1}+1}-\lambda)B}\\
      &\hspace{3cm}
      -\sum_{n_1 > \dots > n_{L_m} > 0}
      \frac{1}{A\,n_{L_{i-1}+1}^{k_i-1}(n_{L_{i-1}+1}-\lambda)B}\\
   \end{aligned}\\
  &\quad=
   \begin{aligned}[t]
    &\sum_{n_1 > \dots > n_{L_m} > 0}
    \frac{1}{A}
    \left\{
    \frac{\lambda}{
                    n_{L_{i-1}+1}^{k_i}(n_{L_{i-1}+1}-\lambda)}
    -\frac{1}{
                 n_{L_{i-1}+1}^{k_i-1}(n_{L_{i-1}+1}-\lambda)}
    \right\} \frac{1}{B}\\
   \end{aligned}\\
  &\quad=
   \sum_{n_1 > \dots > n_{L_m} > 0}
   \frac{1}{A}
    \Biggl\{
     \frac{\lambda'}{(n_{L_{i-1}+1}-1)^{k_i}(n_{L_{i-1}+1}-\lambda)}\\
   &\qquad\left.
      -\frac{1}{(n_{L_{i-1}+1}-1)^{k_i-1}(n_{L_{i-1}+1}-\lambda)}
    -\left(
     \frac{1}{(n_{L_{i-1}+1}-1)^{k_i}} - \frac{1}{n_{L_{i-1}+1}^{k_i}}
    \right)
   \right\} \frac{1}{B}\\
  &\quad=
  \begin{aligned}[t]
   &\lambda' [\dots,l_{i-1},(n-1)^{k_i},l_i,\dots;\lambda]
    - [\dots,l_{i-1},(n-1)^{k_i-1},l_i,\dots;\lambda]\\
   &\quad
   + 
   \!\!\!
   \sum_{\substack{n_1> \dots > n_{L_{i-1}}\\
            n_{L_{i-1}} -1 >n_{L_{i-1}+2} > \dots > n_{L_m} > 0}}
   \!\!\!
   \frac{1}{A}
   \left(
    \frac{1}{(n_{L_{i-1}+2})^{k_i}}
     - \frac{1}{(n_{L_{i-1}}-1)^{k_i}}
   \right) \frac{1}{B}.
  \end{aligned}
 \end{align*}
We divide the range of sum of the third term into two parts as
  \begin{align*}
   \sum_{\substack{n_1> \dots > n_{L_{i-1}}\\
    n_{L_{i-1}} -1 >n_{L_{i-1}+2} > \dots > n_{L_m} > 0}}
   = 
   \sum_{n_1> \dots > n_{L_m} > 0}
    - \sum_{\substack{n_1> \dots > n_{L_{i-1}}\\
       n_{L_{i-1}+2} = n_{L_{i-1}} -1 \\
       n_{L_{i-1}+2} > \dots > n_{L_m} > 0}}.
  \end{align*}
The later sum is equal to zero because of $n_{L_{i-1}+2} = n_{L_{i-1}} -1$.
Thus we have
 \begin{align*}
  &\lambda [\{\dots,l_{i-1},n^{k_i},l_i,\dots\};\lambda]
   - [\{\dots,l_{i-1},n^{k_i-1},l_i,\dots\};\lambda]\\
  &=
  \begin{aligned}[t]
   &\lambda' [\{\dots,l_{i-1},(n-1)^{k_i},l_i,\dots\};\lambda]
    - [\{\dots,l_{i-1},(n-1)^{k_i-1},l_i,\dots\};\lambda]\\
   & + \sum_{n_1> \dots > n_{L_m} > 0}
   \frac{1}{A}
   \left\{
    \frac{1}{(n_{L_{i-1}+2})^{k_i}}
     - \frac{1}{(n_{L_{i-1}}-1)^{k_i}}
   \right\} \frac{1}{B}
  \end{aligned}\\
  &=
  \begin{aligned}[t]
    \lambda'& [\{\dots,l_{i-1},(n-1)^{k_i},l_i,\dots\};\lambda]
    - [\{\dots,l_{i-1},(n-1)^{k_i-1},l_i,\dots\};\lambda]\\
    &- [\{\dots,l_{i-1}-1,(n-1)^{k_i},l_i,\dots\};\lambda]
    +[\{\dots,l_{i-1},n^{k_i},l_i-1,\dots\};\lambda].
  \end{aligned}
 \end{align*}
\item
\begin{enumerate}[{\normalfont \upshape (a)}]
\setcounter{enumii}{2}
\item
 Repeating shift of $n_i \mapsto n_i+1$, we have
 \begin{align*}
  &\begin{aligned}[t]
   &\lambda\,[\{(n-1)^{k_1},l_1,
      \dots,(n-1)^{k_{m-1}},l_{m-1},n^{k_m},1\};\lambda]\\
   &\qquad
   -\lambda\,[\{(n-1)^{k_1},l_1,
      \dots,(n-1)^{k_{m-1}},l_{m-1},n^{k_m-1},1\};\lambda]
  \end{aligned}\\
  &\qquad=
  -\sum_{n_1>\dots>n_{L_m}>0}
   \frac{1}{(n_1-1)^{k_1}\dotsm (n_{L_m-1}-\lambda)n_{L_m}^{k_m}}\\
  &\qquad=
  -\sum_{n_1>\dots>n_{L_m}\ge 0}
   \frac{1}{n_1^{k_1}\dotsm (n_{L_m-1}-\lambda')(n_{L_m}+1)^{k_m}}\\
  &\qquad\qquad\text{(by shift $n_i\longmapsto n_i+1$)}\\
  &\qquad=
  -\sum_{n_1>\dots>n_{L_m-1}\ge n_{L_m}> 0}
   \frac{1}{n_1^{k_1}\dotsm (n_{L_m-1}-\lambda')n_{L_m}^{k_m}}\\
  &\qquad\qquad\text{(by shift $n_L+1 \longmapsto n_L$)}\\
  &\qquad=
  -\sum_{n_1>\dots>n_{L_m}> 0}
   \frac{n_{L_m}-\lambda'}
    {n_1^{k_1}\dotsm(n_{L_m-1}-\lambda')
      n_{L_m}^{k_m}(n_{L_m}-\lambda')}\\
  &\qquad\qquad\qquad\qquad
   -\sum_{n_1>\dots>n_{L_m-1}> 0}
   \frac{1}{n_1^{k_1}\dotsm n_{L_m-1}^{k_m}(n_{L_m-1}-\lambda')}\\
   \begin{split}
    &\qquad=
    \lambda'\,[\{n^{k_1},l_1,\dots,l_{m-1}-1,n^{k_m},1\};\lambda'] \\
    &\qquad\qquad\qquad
     - [\{n^{k_1},l_1,\dots,l_{m-1},n^{k_m-1},1\};\lambda']\\
    &\qquad\qquad\qquad\qquad\qquad\qquad
     - [\{n^{k_1},l_1,\dots,l_{m-1}-1,n^{k_m},1\};\lambda'].
   \end{split} 
 \end{align*} 
\item
Similarly as in the previous cases,
 \begin{align*}
  &[\{(n-1)^{k_i},l_i\}_{i=1}^m;\lambda]\\
  &\quad=
   \sum_{n_1>\dots>n_{L_m}>0}
    \frac{1}{(n_1-1)^{k_1}(n_1-\lambda)
             \dots (n_{L_m-1}-\lambda)(n_{L_m}-\lambda)}\\
  &\quad=
   \sum_{n_1>\dots>n_{L_m}\ge 0}
    \frac{1}{n_1^{k_1}(n_1-\lambda')
             \dots (n_{L_m-1}-\lambda')(n_{L_m}-\lambda')}\\
  &\quad=
%   (\text{$n_{L_m}>0$ $B$N;~(B}) + (\text{$n_{L_m}=0$ $B$N;~(B})\\
   (\text{$n_{L_m}>0$ part}) + (\text{$n_{L_m}=0$ part})\\
  &\quad=
   [\{n^{k_i},l_i\}_{i=1}^m;\lambda']
    - \frac{1}{\lambda'}[\{n^{k_i},l_i\}_{i=1}^{m-1}\cup\{n^{k_m},l_m-1\};\lambda'].
 \end{align*}
\end{enumerate}
\end{enumerate}
\end{proof}

\begin{proof}[Proof of Theorem~\ref{mainTh}]
\begin{enumerate}[{\normalfont \upshape (i)}]
\item
Applying Lemma~\ref{lemma} succesively
 \begin{align*}
 &\text{(LHS)}\\
 &=
  \sum_{\{(\delta_i,\epsilon_i)\}_{i=1}^m\in I^m}
   (-\lambda)^{m-|\delta|-|\epsilon|}\,
   [\{n^{k_i-\delta_i},l_i-\epsilon_i\}_{i=1}^m;\lambda]\\
\begin{split}
 &=
   \sum_{\{(\delta_i,\epsilon_i)\}_{i=2}^m\in I^{m-1}}
  \sum_{\delta_1'\in\{0,1\}}
    (-\lambda)^{m-1-|\delta|-|\epsilon|}(-\lambda')^{1-|\delta'|}\\
    &\qquad\qquad\qquad\times
     [\{(n-1)^{k_1-\delta_1'},l_1,n^{k_2-\delta_2},l_2-\epsilon_2,
      \dots,
      n^{k_m-\delta_m},l_m-\epsilon_m\};\lambda]\\
  &\qquad\qquad\text{(by Lemma~\ref{lemma}~(\ref{lemma.2}))}
\end{split} \\
&=
  \sum_{\{(\delta_i,\epsilon_i)\}_{i=3}^m\in I^{m-2}}
  \sum_{\delta_1'\in\{0,1\}}
  \sum_{(\delta_2',\epsilon_2')\in I}
    (-\lambda)^{m-2-|\delta|-|\epsilon|}
       (-\lambda')^{2-|\delta'|-|\epsilon'|}\\
    &\qquad\qquad\times
    [\{(n-1)^{k_1-\delta_1'},l_1-\epsilon_2',(n-1)^{k_2-\delta_2'},l_2,
    \dots,
    n^{k_m-\delta_m},l_m-\epsilon_m\};\lambda]\\
  &\qquad\qquad\text{(by Lemma~\ref{lemma}~(\ref{lemma.1}))}\\
&=
  \sum_{\delta_1'\in\{0,1\}}
  \sum_{\{(\delta_j',\epsilon_j')\}_{j=2}^m\in I^{m-1}}
    (-\lambda')^{m-|\delta'|-|\epsilon'|}\\
    &\quad\quad\times
    [\{(n-1)^{k_1-\delta_1'},l_1-\epsilon_2',
     (n-1)^{k_2-\delta_2'},l_2-\delta_3',
      \dots,
       (n-1)^{k_m-\delta_m'},l_m\};\lambda] \\
  &\qquad\qquad\text{(by Lemma~\ref{lemma}~(\ref{lemma.1}) $m-2$ times)}\\
&=
  \sum_{\delta_1',\epsilon_{m+1}'\in\{0,1\}}
  \sum_{\{(\delta_j',\epsilon_j')\}_{j=2}^m\in I^{m-1}}
    (-\lambda')^{m-|\delta'|-|\epsilon'|}\\
    &\qquad\qquad\qquad\times
     [\{n^{k_1-\delta_1'},l_1-\epsilon_2',n^{k_2-\delta_2'},l_2-\delta_3',
      \dots,
       n^{k_m-\delta_m'},l_m-\epsilon_{m+1}\};\lambda']\\
  &\qquad\qquad\text{(by Lemma~\ref{lemma}~(\ref{lemma.5}))}\\
&=
  \text{(RHS).}
 \end{align*} 
\item
 Use 
  Lemma~\ref{lemma}~(\ref{lemma.3}),
  Lemma~\ref{lemma}~(\ref{lemma.1}) $m-1$ times,
  and Lemma~\ref{lemma}~(\ref{lemma.5}).
\end{enumerate}

\vspace{0.5cm}
Remaining relations
  and the relations of $g$'s can be proved quite similarly.
\end{proof}

\section{Another proof of Ohno's relation.\label{sec.proof}}
 First of all, we introduce a partial order of indices
 $\{k_1,l_1,\dotsc,k_m,l_m\}\in\Z_{\ge 0}^{2m}$ as follows:
 \begin{align*}
  &\{k_1,l_1,k_2,l_2,\dots,k_m,l_m\}
   \le \{k_1',l_1',k_2',l_2',\dots,k_{m'}',l_{m'}'\}\\
  &\qquad\stackrel{\mathrm{def}}{\Longleftrightarrow}
    \quad\left\{\begin{aligned}[c]
     &m < m' \qquad \text{ or }\\
     &m = m \quad\text{ and }\quad k_i'-k_i,\ l_i'-l_i \ge 0 \quad
	  \text{for all $i$}.
     \end{aligned}   
    \right.
 \end{align*}
 
By induction with respect to this order, we prove Ohno's relation.
 If the index is minimum i.e. $(1,1)$, it is clear because of $\tau(xy) = xy$.

 If the theorem is correct less than $\{k_i,l_i\}_{i=1}^m$,
 applying Theorem~\ref{mainTh} to
 $f(\{k_i,l_i\}_{i=1}^m;\lambda)$ and $g(\{k_i,l_i\}_{i=1}^m;\lambda)$,
 we obtain two relations for $f$ and $g$.
 Subtracting these two equations, we have
 \begin{align*}
  &\sum{\lambda}^{m-|\delta| - |\epsilon|}\,
   \Bigl\{
   f(\{k_i-\delta_i,l_i-\epsilon_i\}_{i=1}^m;\lambda)
   - g(\{k_i-\delta_i,l_i-\epsilon_i\}_{i=1}^m;\lambda)
  \Bigr\}
  \\
  &\quad=
  \sum{\lambda'}^{m-|\delta'|-|\epsilon'|}
   \Bigl\{
    f(\{k_i-\delta_i',l_i-\epsilon_i'\}_{i=1}^m;\lambda')
   - g(\{k_i-\delta_i',l_i-\epsilon_i'\}_{i=1}^m;\lambda')
   \Bigr\}.
 \end{align*}
 But the terms whose indices are less than $\{k_i,l_i\}_{i=1}^m$
 are canceled out
 by the induction hypothesis.
 The remaining is
 \begin{equation*}
 \begin{split}
  &{\lambda}^m
   \Bigl\{
    f(\{k_i,l_i\}_{i=1}^m;\lambda)
   - g(\{k_i,l_i\}_{i=1}^m;\lambda)
   \Bigr\}
   \\
  &\qquad=
  {\lambda'}^m
  \Bigl\{
   f(\{k_i,l_i\}_{i=1}^m;\lambda')
   - g(\{k_i,l_i\}_{i=1}^m;\lambda')
  \Bigr\}.
 \end{split} 
 \end{equation*}
 Hence
  ${\lambda}^m\,f(\{k_i,l_i\}_{i=1}^m;\lambda)
    - {\lambda}^m\,g(\{k_i,l_i\}_{i=1}^m;\lambda) $
 is a periodic function in $\lambda$ with period $1$.
 Furthermore by Lemma~\ref{key.lem} it has 
 partial-fractions expansion such as
 \begin{align*}
  \lambda^m \sum_{n=1}^\infty
   \frac{C_n}{n-\lambda}.
 \end{align*}
 So it must be zero.
 Thus we complete the proof.

\section*{Acknowledgements}

 The authors express their deep gratitude
 to Professor Masanobu Kaneko and Professor Shigeki Akiyama
 for their valuable suggestions.

 The second author is partially supported by
 Grant-in-Aid Scientific Research from the Ministry of 
 Education, Culture, Sports, Science and Technology
 (12640046) and by Waseda University Grant for Special Research 
 Project (2000A-124).

\begin{tabbing}
 Department of Mathematical Sciences\\
 School of Science and Engineering,\\
 WASEDA UNIVERSITY\\
 3-4-1, Okubo Shinjuku-ku\\
 Tokyo 169-8555, Japan\\
  E-mail: \= okuda@gm.math.waseda.ac.jp\\
          \> uenoki@mn.waseda.ac.jp
\end{tabbing}

\end{document}